\documentclass[11pt,a4paper,leqno,twoside]{amsart}
\usepackage{latexsym,amssymb,amsmath}
\input xy
\xyoption{all}

\def\CC{\mathbb C}

\def\RR{\mathbb R}
\def\HH{\mathbb H}
\def\AA{{\mathbb A}}

\def\OO{\mathbb O}

\def\11{\mathbf 1}
\def\PP{\mathbb P}

\def\e1{\varepsilon_1}
\def\e2{\varepsilon_2}
\def\e3{\varepsilon_3}

\def\P2{{\PP}^2}

\def\00{\underline{0}}
\def\J0{{\cal J}_3(\underline{0})}

\def\PJ0{\PP({\cal J}_3(\underline{0}))}

\def\e{\varepsilon}

\def\AP2{{\AA\PP}^2}
\def\RP2{{\RR\PP}^2}
\def\CP2{{\CC\PP}^2}
\def\HP2{{\HH\PP}^2}
\def\OP2{{\OO\PP}^2}

\newtheorem{theo}{Theorem}[section]
\newtheorem{coro}[theo]{Corollary}
\newtheorem{lemm}[theo]{Lemma}
\newtheorem{prop}[theo]{Proposition}

\newtheorem{prob}[theo]{Problem}

\theoremstyle{remark}
\newtheorem{rema}[theo]{Remark}

\begin{document}
\title[Fields of dimension one]{Fields of dimension one algebraic
over a global or local field need not be of type $C _{1}$}
\keywords{Field of dimension $\le 1$, field of type $C_{1}$, form,
Henselian valuation\\
2020 MSC Classification: 11E76 11R34 (primary), 12F10, 12J10, 11S15
(secondary).}

\author{Ivan D. Chipchakov}
\address{Institute of Mathematics and Informatics\\Bulgarian Academy
of Sciences\\1113 Sofia, Bulgaria: E-mail address:
chipchak@math.bas.bg}

\begin{abstract}
Let $(K, v)$ be a Henselian discrete valued field with a quasifinite
residue field. This paper proves the existence of an algebraic
extension $E/K$ satisfying the following: (i) $E$ has dimension
dim$(E) \le 1$, i.e. the Brauer group Br$(E ^{\prime })$ is trivial,
for every algebraic extension $E ^{\prime }/E$; (ii) finite
extensions of $E$ are not $C _{1}$-fields. This, applied to the
maximal algebraic extension $K$ of the field $\mathbb{Q}$ of rational
numbers in the field $\mathbb{Q} _{p}$ of $p$-adic numbers, for a
given prime $p$, proves the existence of an algebraic extension $E
_{p}/\mathbb{Q}$, such that dim$(E _{p}) \le 1$, $E _{p}$ is not a $C
_{1}$-field, and $E _{p}$ has a Henselian valuation of residual
characteristic $p$.
\end{abstract}

\maketitle

\par
\medskip
\section{\bf Introduction}

\medskip
A field $F$ is said to be of dimension $\le 1$, if the Brauer groups
Br$(F ^{\prime })$ are trivial, for all algebraic field extensions $F
^{\prime }/F$. It is known (cf. \cite{S1}, Ch.~II, 3.1) that dim$(F)
\le 1$ if and only if Br$(F ^{\prime }) = \{0\}$, where $F ^{\prime }$
runs across the set Fe$(F)$ of finite extensions of $F$ in its
separable closure $F _{\rm sep}$. When $F$ is perfect, we have
dim$(F) \le 1$ if and only if the absolute Galois group
$\mathcal{G}_{F}$ $:= \mathcal{G}(F _{\rm sep}/F)$ has cohomological
dimension cd$(\mathcal{G}_{F}) \le 1$ as a profinite group. For
example, we have dim$(F) \le 1$ in case $F$ is a quasifinite field,
i.e. a perfect field which admits, for each $n \in \mathbb{N}$, a
unique extension in $F _{\rm sep}$ of degree $n$. It is well-known
that then $\mathcal{G}_{F}$ is isomorphic to
$\mathcal{G}_{\mathbb{F}}$, for any finite field $\mathbb{F}$, i.e.
it is a profinite completion of the additive group $\mathbb{Z}$ of
integers. One has cd$(\mathcal{G}_{F}) \le 1$, since
$\mathcal{G}_{F}$ is isomorphic to the topological group product
$\prod _{p \in \mathbb{P}} \mathbb{Z} _{p}$, where $\mathbb{P}$ is
the set of prime numbers and $\mathbb{Z} _{p}$ is the additive group
of $p$-adic integers, for each $p$ (see \cite{GiSz}, Examples~4.1.2).
\par
We say that $F$ is of type $C _{r}$ (or a $C _{r}$-field), for some
$r \in \mathbb{N}$, if  every $F$-form (a homogeneous nonzero
polynomial with coefficients in $F$) $f$ of degree deg$(f)$ in more
than deg$(f) ^{r}$ variables has a nontrivial zero over $F$. For $F$
to be a $C _{r}$-field, it is necessary that char$(F) = q$, and
in case $q > 0$, the degree $[F\colon F ^{q}]$ of $F$ as an
extension of its subfield $F ^{q} = \{\beta ^{q}\colon \beta \in F\}$
be at most equal to $q ^{r}$. The class of $C _{r}$-fields is closed
under the formation of algebraic extensions, and it contains the
extensions of transcendency degree $r$ over any algebraically closed
field (cf. \cite{L1}). It is known that $C _{1}$-fields have
dimension $\le 1$ but the converse is not necessarily true in nonzero
characteristic; \label{k999} one can take as a counter-example any
field $\Phi $ with char$(\Phi ) = q > 0$, $[\Phi \colon \Phi ^{q}] >
q$ and $\Phi _{\rm sep} = \Phi $ (see \cite{S1}, Ch. II, 3.1 and 3.2).
The restriction on the characteristic has been lifted by Ax \cite{Ax1},
by providing an example of a quasifinite field $F$ (so having dim$(F)
\le 1$), such that char$(F) = 0$ and $F$ is not a $C _{r}$-field, for any $r \in
\mathbb{N}$. This example is in sharp contrast to the
Chevalley-Warning theorem (see, e.g., \cite{GiSz}, Theorem~6.2.6),
which establishes the $C _{1}$ type of finite fields.
\par
\medskip
As noted by Serre in \cite{S1}, Ch. II, 3.3, it is not known whether
an algebraic extension $E$ of the field $\mathbb{Q}$ of rational
numbers is of type $C _{1}$, provided that dim$(E) \le 1$; he has
added that this is not likely to hold in general.
\par
\label{k999}
The present paper answers the question arising from Serre's remark,
in the direction pointed there. The answer is unchanged when
$\mathbb{Q}$ is replaced by any global or local field. This is
obtained by methods of valuation theory, using the fact that
nontrivial Krull valuations of global fields are discrete with finite
residue fields (see \cite{E3}, Examples~4.1.2, 4.1.3 and
Corollary~14.2.2). The considered question remains open for Galois
extensions $E$ of $\mathbb{Q}$ with dim$(E) \le 1$.
\par
\medskip
\section{\bf Statements of the main results}
\par
\medskip
The main results of this paper are presented as two theorems. The
former theorem is a special case of the latter one and can be stated
as follows:
\par
\medskip
\begin{theo}
\label{theo2.1}
For each prime number $p$, there exists an algebraic extension $E
_{p}$ of the field $\mathbb{Q}$ of rational numbers, such that {\rm
dim}$(E _{p}) \le 1$, the finite extensions of $E _{p}$ are not $C
_{1}$-fields, and $E _{p}$ is endowed with a Henselian valuation
$v$ whose residue field $\widehat E _{p}$ is of characteristic $p$.
\end{theo}
\par
\medskip
As usual, by a Henselian valuation of a field $K$, we mean a
nontrivial Krull valuation $v$ that extends uniquely, up-to
equivalence, to a valuation $v _{L}$ on each algebraic extension $L$
of $K$. When $v$ is Henselian, $(K, v)$ is called a Henselian field.
Our next result, stated below, shows that the class of separable
(algebraic) extensions of Henselian discrete valued\footnote{In what
follows, we write briefly "HDV" instead of "Henselian discrete
valued".} fields with finite residue fields contains fields of
dimension $\le 1$, which are not of type $C _{1}$.
\par
\medskip
\begin{theo}
\label{theo2.2} Let $(K, v)$ be an {\rm HDV}-field with a quasifinite
residue field. Then there exists an extension $E$ of $K$ in $K _{\rm
sep}$, such that {\rm dim}$(E) \le 1$ and finite extensions of $E$
are not $C _{1}$-fields. Moreover, $E$ can be chosen so that there is
a sequence $f _{n}$, $n \in \mathbb{N}$, of $E$-forms subject to the
following restrictions:
\par
{\rm (a)} $f _{n}$ does not possess a nontrivial zero over $E$, for
any index $n$;
\par
{\rm (b)} The degrees {\rm deg}$(f _{n}) := p _{n}$, $n \in \mathbb
N$, form a strictly increasing sequence of prime numbers, and for
each $n$, $f _{n}$ depends essentially on exactly $p _{n}k_{n}$
variables, for some $k _{n} \in \mathbb{N}$ with $2 \le k _{n} \le (p
_{n} - 1)/2$.
\end{theo}
\par
\medskip
Theorem \ref{theo2.2} is proved in Section 4. Here we show that
Theorem \ref{theo2.2} implies Theorem \ref{theo2.1}. Denote by
$\mathbb{P}$ the set of prime numbers, and for each $p \in
\mathbb{P}$, let $K _{p}$ be the maximal separable extension of
$\mathbb{Q}$ in the field $\mathbb{Q} _{p}$ of $p$-adic numbers.
It is known (cf. \cite{E3}, Theorem~15.3.5) that the valuation,
say $\omega _{p}$, induced on $K _{p}$ by the standard valuation
of $\mathbb{Q} _{p}$ is Henselian and discrete, and the residue
field $\widehat K _{p}$ of $(K _{p}, \omega _{p})$ is a field with
$p$ elements. Hence, by Theorem \ref{theo2.2}, it has an algebraic
extension $E _{p}$ with the properties required by Theorem
\ref{theo2.1}. The question of whether an algebraic extension $E$
of $\mathbb{Q}$ with dim$(E) \le 1$ is a $C _{\nu }$-field, for
some integer $\nu \ge 2$, remains open. Note in this connection
that examples given by Arkhipov and Karatsuba \cite{AK} (see also
\cite{BS}, Ch. I, Sect. 6.5, and further references there) show
that $\mathbb{Q} _{p}$, $p \in \mathbb{P}$, are not $C _{\nu
}$-fields, for any $\nu \in \mathbb{N}$.
\par
Similarly, if $\mathbb{F} _{p}(X)$ is the rational function field in
a variable $X$ over the field $\mathbb{F} _{p}$ with $p$ elements, $K
_{p}$ is the maximal separable extension of $\mathbb{F} _{p}(X)$ in
the formal Laurent power series field $\mathbb{F} _{p}((X))$, and
$\omega _{p}$ is the valuation of $K _{p}$ induced by the natural
discrete valuation of $\mathbb{F} _{p}((X))$, then $(K _{p}, \omega
_{p})$ is an HDV-field. Therefore, any extension $E _{p}$ of $K _{p}$
in $K _{p,{\rm sep}}$ with the properties claimed by Theorem
\ref{theo2.2} is an algebraic extension of $\mathbb{F} _{p}(X)$, such
that dim$(E _{p}) \le 1$ and finite extensions of $E _{p}$ are not $C
_{1}$-fields ($E _{p}$ is a $C _{2}$-field, by \cite{L1}, Theorem~8).
Note here that $[E _{p}\colon E _{p} ^{p}] = p$ (see (3.2) (a) and
(3.5) (b)), i.e. $E _{p}$ satisfies the necessary condition for
having type $C _{1}$, stated in Section 1 (and proved in \cite{S1},
Ch. II, 3.2). These facts (see also (3.5) (a) and the proof of (3.5))
show that $C _{1}$-fields, global fields and local fields are almost
perfect, in the following sense:
\par
\medskip\noindent
A field $E$ is called almost perfect if one of the following two
equivalent conditions holds: (i) every finite extension of $E$ has a
primitive element; (ii) char$(E) = q$, and in case $q > 0$, $[E\colon
E ^{q}]$ equals $1$ or $q$ (the equivalence of conditions (i) and
(ii) is implied by \cite{L}, Ch. V, Theorem~4.6 and Corollary~6.10).
\par
\medskip
Theorem \ref{theo2.2} and the method of proving Theorem
\ref{theo2.1} enable one to answer the main question considered in
this paper, for Henselian fields that are algebraic extensions of
an arbitrary global field $K$ (see Corollary \ref{coro4.2}). At
the same time, the discussed question remains widely open over
other interesting classes of algebraic extensions of $K$, e.g.,
the class $\mathcal{V}(K)$ of algebraic extensions $E/K$, such
that $\mathcal{G}_{E}$ is a torsion-free group and Henselizations
of $E$ in $E _{\rm sep}$ with respect to any nontrivial valuation
coincide with $E _{\rm sep}$. By the Artin-Schreier theory (cf.
\cite{L}, Ch. XI, Theorem~2.2), the condition on $\mathcal{G}_{E}$
restates the one that $E$ is a nonreal field, i.e. $-1$ is
presentable over $E$ as a finite sum of squares; its violation
yields Br$(E) \neq \{0\}$ (the equivalence class in Br$(E)$ of the
$E$-algebra of Hamiltonian quaternions is an element of order
$2$). When $E$ is nonreal, the description of Br$(K)$ by class
field theory (cf. \cite{We}, Ch. XIII, Sects. 3, 6) enables one to
obtain that dim$(E) \le 1$ by the method of proving Proposition~9
of \cite{S1}, Ch. II. Also, by the Frey-Prestel theorem (see
\cite{Fr}, Theorem~2, and \cite{FJ}, Corollary~11.5.5),
PAC\footnote{As usual, PAC is an abbreviation for "pseudo
algebraically closed".} fields are contained in the class
$\mathcal{V}$ of nonreal fields whose Henselizations with respect
to any nontrivial valuation are separably closed. The main result
of \cite{GJ} shows the existence of non-PAC fields in
$\mathcal{V}$ and raises interest in the following open problem
(see Remark \ref{rema4.4} and \cite{FJ}, Problem~11.5.9~(b)):
\par
\medskip
\begin{prob}
\label{prob}
For a global field $K$, find whether all $\Phi \in \mathcal{V}(K)$
are PAC fields.
\end{prob}
\par
\medskip
The basic notation, terminology and conventions kept in this paper
are standard and virtually the same as in \cite{TW}, \cite{L} and
\cite{Ch1}. Throughout, Brauer and value groups are written
additively, Galois groups are viewed as profinite with respect to
the Krull topology, and by a profinite group homomorphism, we mean
a continuous one.  For any field $E$, $E ^{\ast }$ is its
multiplicative group, $E ^{\ast n} = \{a ^{n}\colon a \in E ^{\ast
}\}$, for each $n \in \mathbb{N}$, and Br$(E) _{p}$, $p \in
\mathbb{P}$, are the $p$-components of Br$(E)$. As usual, for any
$p \in \mathbb{P}$, $E(p)$ denotes the maximal $p$-extension of
$E$ (in $E _{\rm sep}$), that is, the compositum of those finite
Galois extensions of $E$ in $E _{\rm sep}$, whose Galois groups
are $p$-groups. Given a field extension $E ^{\prime }/E$, we write
I$(E ^{\prime }/E)$ for the set of intermediate fields of $E
^{\prime }/E$, $\pi _{E\to E'}$ for the scalar extension map
Br$(E) \to {\rm Br}(E ^{\prime })$, and Br$(E ^{\prime }/E)$ for
the relative Brauer group of $E ^{\prime }/E$ (the kernel of $\pi
_{E\to E'}$). When $E ^{\prime }/E$ is a finite extension, $N
_{E'/E}$ denotes the norm map $E ^{\prime } \to E$, and $N(E
^{\prime }/E)$ stands for the norm group of $E ^{\prime }/E$ (the
image of $E ^{\prime \ast }$ under $N _{E'/E}$). Moreover, if $E
^{\prime }/E$ is a Galois extension, then its Galois group is
denoted by $\mathcal{G}(E ^{\prime }/E)$; we say that $E ^{\prime
}/E$ is a cyclic extension if $\mathcal{G}(E ^{\prime }/E)$ is a
cyclic group. By a $\mathbb{Z} _{p}$-extension, we mean a Galois
extension $\Psi ^{\prime }/\Psi $ with $\mathcal{G}(\Psi ^{\prime
}/\Psi )$ isomorphic to $\mathbb{Z} _{p}$. The value group of any
discrete valued field is assumed to be an ordered subgroup of the
additive group of the field $\mathbb{Q}$; this is done without
loss of generality, in view of \cite{E3}, Theorem~15.3.5, and the
fact that $\mathbb{Q}$ is a divisible hull of any of its infinite
subgroups (see page \pageref{k99}).
\par
\medskip
Here is an overview of this paper: Section 3 includes
valuation-theoretic preliminaries used in the sequel as well as
characterizations of fields of dimension $\le 1$ among algebraic
extensions of local fields. As noted above, Theorem \ref{theo2.2} is
proved in Section 4. This is done by modifying the proof of the
Theorem of \cite{Ax1}, given in \cite{Ax2}; specifically, the forms
violating the $C _{1}$ condition are defined by essentially the same
pattern in both proofs.
\par
\medskip
\section{\bf Preliminaries and characterizations of algebraic
extensions $E$ of local fields with Br$(E) _{p} = \{0\}$, for a
given prime $p$}
\par
\medskip
For any field $K$ with a (nontrivial) Krull valuation $v$, $O
_{v}(K) = \{a \in K\colon \ v(a) \ge 0\}$ denotes the valuation ring
of $(K, v)$, $M _{v}(K) = \{\mu \in K\colon \ v(\mu ) > 0\}$ the
maximal ideal of $O _{v}(K)$, $O _{v}(K) ^{\ast } = \{u \in K\colon
\ v(u) = 0\}$ the multiplicative group of $O _{v}(K)$, $v(K)$ the
value group and $\widehat K = O _{v}(K)/M _{v}(K)$ the residue field
of $(K, v)$, respectively; $\overline {v(K)}$ is a divisible hull of
$v(K)$. The condition that $v$ is Henselian has the following two
equivalent forms (cf. \cite{E3}, Sect. 18.1):
\par
\medskip\noindent
(3.1) (a) Given a polynomial $f(X) \in O _{v}(K) [X]$ and an element
$a \in O _{v}(K)$, such that $2v(f ^{\prime }(a)) < v(f(a))$, where
$f ^{\prime }$ is the formal derivative of $f$, there is a zero $c
\in O _{v}(K)$ of $f$ satisfying the equality $v(c - a) = v(f(a)/f
^{\prime }(a))$;
\par
(b) For each normal extension $\Omega /K$, $v ^{\prime }(\tau (\mu ))
= v ^{\prime }(\mu )$ whenever  $\mu \in \Omega $, $v ^{\prime }$ is
a valuation of $\Omega $ extending $v$, and $\tau $ is a
$K$-automorphism of $\Omega $.
\par
\medskip
Next we recall some facts concerning the case where $(K, v)$ is a
real-valued field, i.e. $v(K)$ is embeddable as an ordered subgroup
in the additive group $\mathbb{R}$ of real numbers. Fix a completion
$K _{v}$ of $K$ with respect to the topology induced by $v$, and
denote by $\bar v$ the valuation of $K _{v}$ continuously extending
$v$. Then:
\par
\medskip\noindent
(3.2) (a) $(K, v)$ is Henselian if and only if $K$ has no proper
separable (algebraic) extension in $K _{v}$ (cf. \cite{E3},
Corollary~18.3.3);
\par
(b) The topology of $K _{v}$ as a completion of $K$ is the same as
the one induced by $\bar v$; also, $\bar v(K _{v}) = v(K)$, $\widehat
K$ equals the residue field of $(K _{v}, \bar v)$, and $(K _{v}, \bar
v)$ is a Henselian field (cf. \cite{E3}, Theorems~9.3.2 and 18.3.1).
\par
\medskip
When $v$ is Henselian (but not necessarily real-valued), so is $v
_{L}$, for any algebraic field extension $L/K$. In this case, we
denote by $\widehat L$ the residue field of $(L, v _{L})$, and put
$O _{v}(L) = O _{v _{L}}(L)$, $M _{v}(L) = M _{v_{L}}(L)$, $v(L) = v
_{L}(L)$; also, we write $v$ instead of $v _{L}$ when there is no
danger of ambiguity. Clearly, $\widehat L/\widehat K$ is an algebraic
extension and $v(K)$ is an ordered subgroup of $v(L)$, such that
$v(L)/v(K)$ is a torsion group; hence, one may assume without loss of
generality that \label{k99} $v(L)$ is an ordered subgroup of
$\overline {v(K)}$. By Ostrowski's theorem (cf. \cite{E3},
Theorem~17.2.1), if $[L\colon K]$ is finite, then it is divisible by
$[\widehat L\colon \widehat K]e(L/K)$, and in case $[L\colon K] \neq
[\widehat L\colon \widehat K]e(L/K)$, the integer $[L\colon K][\widehat
L\colon \widehat K] ^{-1}e(L/K) ^{-1}$ is a power of char$(\widehat
K)$ (so char$(\widehat K) \mid [L\colon K]$); here $e(L/K)$ is the
ramification index of $L/K$, i.e. the index of $v(K)$ in $v(L)$. Ostrowski's
theorem implies the following:
\par
\medskip\noindent
(3.3) The quotient groups $v(K)/pv(K)$ and $v(L)/pv(L)$ are
isomorphic, if $p \in \mathbb{P}$ and $[L\colon K] < \infty $.
Moreover, if char$(\widehat K) \nmid [L\colon K]$, then the natural
embedding of $K$ into $L$ induces canonically an isomorphism
$v(K)/pv(K) \cong v(L)/pv(L)$.
\par
\medskip
The finite extension $L/K$ satisfies the equality $[L\colon K] =
[\widehat L\colon \widehat K]e(L/K)$ in each of the following two
situations:
\par
\medskip\noindent
(3.4) (a) $(K, v)$ is HDV and $L/K$ is separable (see \cite{E3},
Sect. 17.4).
\par
(b) $(K, v)$ is HDV and the field $K$ is almost perfect (cf.
\cite{L}, Ch. XII, Proposition~6.1). When char$(K) = q > 0$, this
implies $[K\colon K ^{q}] = q$, since $K ^{q}$ does not contain any
$\alpha \in K ^{\ast }$ with $v(\alpha ) \notin qv(K)$.
\par
\medskip
We show that the conditions of (3.4) (b) hold in the following two
cases:
\par
\medskip\noindent
(3.5) (a) $(K, v)$ is a complete discrete valued field (i.e. $(K, v)
= (K _{v}, \bar v)$ and $v$ is discrete) with $\widehat K$ perfect.
\par
(b) $(K, v)$ is HDV and $K$ is an algebraic extension of a global
field $K _{0}$ (see \cite{E3}, Example~4.1.3).
\par
\medskip\noindent
Complete discrete valued fields are Henselian, by (3.2), so it
suffices for the proof of (3.5) to show that $[K\colon K ^{q}] = q$,
provided that char$(K) = q > 0$. Under the hypotheses of (3.5) (a),
$K$ is isomorphic to the formal Laurent power series field $\widehat
K((X))$ in a variable $X$ over $\widehat K$ (see \cite{E3},
Theorem~12.2.3). Since $\widehat K$ is perfect, whence,
$\widehat K((X)) ^{q} = \widehat K((X ^{q}))$, and $1, X, \dots , X
^{q-1}$ is a basis of $\widehat K((X))$ over $\widehat K((X
^{q}))$, this yields $[K\colon K ^{q}] = q$. Similarly, it turns out
that if $K _{0} = \widehat K(X)$ and $\widehat K$ is the prime field
$\mathbb{Z}/q\mathbb{Z}$, then $K _{0} ^{q} = \widehat K(X ^{q})$ and
$1, X, \dots , X ^{q-1}$ is a basis of $K _{0}$ over $K _{0} ^{q}$.
In view of \cite{BH}, Lemma~2.12, this means that $[K _{0}\colon K
_{0} ^{q}] = q$, for any global field $K _{0}$ of characteristic $q$, and
also proves (3.5) (b).
\par
\medskip
Assume now that $(K, v)$ is a Henselian field and let $R$ be a
finite extension of $K$. We say that $R/K$ is inertial, if
$[R\colon K] = [\widehat R\colon \widehat K]$ and $\widehat R$ is
separable over $\widehat K$; $R/K$ is called totally ramified, if
$e(R/K) = [R\colon K]$. Inertial extensions of $K$ are clearly
separable; also, they have a number of useful properties, some of
which are presented by the following lemma (for its proof, see
\cite{TW}, Theorem~A.23):
\par
\medskip
\begin{lemm}
\label{lemm3.1}
Let $(K, v)$ be a Henselian field and $K _{\rm ur}$ the compositum of
inertial extensions of $K$ in $K _{\rm sep}$. Then:
\par
{\rm (a)} An inertial extension $R ^{\prime }/K$ is Galois if and
only if so is $\widehat R ^{\prime }/\widehat K$. When this holds,
$\mathcal{G}(R ^{\prime }/K)$ and $\mathcal{G}(\widehat R ^{\prime
}/\widehat K)$ are canonically isomorphic.
\par
{\rm (b)} $v(K _{\rm ur}) = v(K)$ and $K _{\rm ur}/K$ is a Galois
extension with $\mathcal{G}(K _{\rm ur}/K) \cong
\mathcal{G}_{\widehat K}$.
\par
{\rm (c)} Finite extensions of $K$ in $K _{\rm ur}$ are inertial, and
the natural mapping of $I(K _{\rm ur}/K)$ into $I(\widehat K _{\rm
sep}/\widehat K)$, by the rule $L \to \widehat L$, is bijective.
\end{lemm}
\par
\medskip
The next two lemmas enable one to generalize a number of results
on complete real-valued fields to the case of Henselian real-valued
fields.
\par
\medskip
\begin{lemm}
\label{lemm3.2} Let $(K, v)$ be a real-valued field, $(K _{v}, \bar
v)$ its completion, and $(K ^{\prime }, v')$ an intermediate valued
field of $(K _{v}, \bar v)/(K, v)$. Suppose that $(K ^{\prime }, v')$
is Henselian, identify $K ^{\prime } _{\rm sep}$ with its $K
^{\prime }$-isomorphic copy in $K _{v,{\rm sep}}$, and denote by $f$
the mapping $f: {\rm Fe}(K ^{\prime }) \to {\rm Fe}(K
_{v})$, defined by the rule $\Lambda ^{\prime } \to \Lambda ^{\prime
}K _{v}$. Then:
\par
{\rm (a)} $K ^{\prime } _{\rm sep} \cap K _{v} = K ^{\prime }$, and
each $\Lambda \in {\rm Fe}(K _{v})$ contains a primitive element
$\lambda \in K ^{\prime } _{\rm sep}$ over $K _{v}$, such that $[K
_{v}(\lambda )\colon K _{v}] = [K ^{\prime }(\lambda )\colon K
^{\prime }]$;
\par
{\rm (b)} $K ^{\prime } _{\rm sep}K _{v} = K _{v,{\rm sep}}$ and
$\mathcal{G}_{K'} \cong \mathcal{G}_{K _{v}}$;
\par
{\rm (c)} The mapping $f$ is bijective and degree-preserving.
Moreover, $f$ and the inverse mapping $f ^{-1}: {\rm Fe}(K _{v}) \to
{\rm Fe}(K ^{\prime })$, preserve the Galois property and the
isomorphism class of the corresponding Galois groups.
\end{lemm}
\par
\medskip
\begin{proof}
The conditions on $(K, v)$ and the Henselian property of $(K ^{\prime
}, v')$ ensure that $K ^{\prime } _{\rm sep} \cap K _{v} = K ^{\prime
}$. The latter part of Lemma \ref{lemm3.2} (a) can be deduced from
Krasner's lemma (see \cite{L2}, Ch. II, Propositions~3, 4). The
conclusions of Lemma \ref{lemm3.2} (c) follow from Lemma
\ref{lemm3.2} (a) and Galois theory (cf. \cite{L}, Ch. VI,
Theorem~1.12), and those of Lemma \ref{lemm3.2} (b) follow from Lemma
\ref{lemm3.2} (a), (c) and the definition of the Krull topology on
$\mathcal{G}_{K'}$ and $\mathcal{G}_{K _{v}}$.
\end{proof}
\par
\medskip
\begin{lemm}
\label{lemm3.3} Let $(K, v)$ be a Henselian real-valued field, $(K
_{v}, \bar v)$ its completion, and $R$ an extension of $K$ in $K
_{\rm sep}$. Identify $K _{\rm sep}$ with its $K$-isomorphic copy
in $R _{v,{\rm sep}}$, where $R _{v} = R _{v _{R}}$, and $K _{v}$
with the topological closure of $K$ in $R _{v}$, and put $R
^{\prime } = K _{v}R$. Then $(R ^{\prime }, \bar v _{R'})$ is an
intermediate valued field of $(R _{v}, \overline {v _{R}})/(R,
v_{R})$.
\end{lemm}
\par
\medskip
\begin{proof}
It follows from Lemma \ref{lemm3.2} (a), (c) and the Henselian
property of $(K, v)$ that the mapping Fe$(K) \to {\rm Fe}(K
_{v})$, by the rule $\Lambda \to K _{v}\Lambda $, is bijective and
degree-preserving. This implies that, for each $\Lambda \in {\rm
Fe}(K)$, the restriction of the norm map $N _{K _{v}\Lambda /K
_{v}}$ on $\Lambda $ equals $N _{\Lambda /K}$, which shows that
$(K _{v}\Lambda , \bar v _{K _{v}\Lambda })/(\Lambda , v _{\Lambda
})$ is a valued field extension (see, e.g., \cite{TW}, Lemma~1.6).
At the same time, observing that $(K _{v},\bar v)$ is Henselian
and $K _{v}\Lambda $ is a completion of $\Lambda $ with respect to
the topology of $v _{\Lambda }$ (see \cite{L}, Ch. XII,
Proposition~3.1), one obtains that if $\Lambda $ is a finite
extension of $K$ in $R$, then $\overline {v _{R}}$ extends $\bar v
_{K _{v}\Lambda } = \overline {v _{\Lambda }}$ upon $R _{v}$. As
$R ^{\prime }$ equals the union $\cup K _{v}\Lambda $, when
$\Lambda $ runs across the set of finite extensions of $K$ in $R$,
these facts prove Lemma \ref{lemm3.3}.
\end{proof}
\par
\medskip
Next we present characterizations of those fields of dimension $\le 1$,
which lie in the class of algebraic extensions of any HDV-field $(K,
v)$ with $\widehat K$ quasifinite. They are stated as two lemmas. For
a proof of the former one in the case where $(K, v)$ is a local field
with char$(K) = 0$, see \cite{S1}, Ch. II, 5.6, Lemma~3.
\par
\medskip
\begin{lemm}
\label{lemm3.4} Let $(K, v)$ be an {\rm HDV}-field with $\widehat K$
quasifinite, and let $R/K$ be an algebraic field extension. Fix some
$p \in \mathbb{P}$, and in case $p = {\rm char}(K)$, suppose that
$R/K$ is separable or $K$ is an almost perfect field. Then {\rm
Br}$(R) _{p} = \{0\}$ if and only if one of the following three
equivalent conditions is fulfilled:
\par
{\rm (a)} {\rm Br}$(R ^{\prime }) _{p} = \{0\}$, for every algebraic
extension $R ^{\prime }/R$; this holds if and only if {\rm Br}$(R
^{\prime }) _{p} = \{0\}$, when $R ^{\prime }$ runs across the set
{\rm Fe}$(R)$;
\par
{\rm (b)} For any pair $R _{1} ^{\prime } \in {\rm Fe}(R)$, $R
^{\prime } \in I(R _{1} ^{\prime }/R)$, $p$ does not divide the
period of the quotient group of $R ^{\prime \ast }$ by the norm group
$N(R _{1} ^{\prime }/R ^{\prime })$ of the extension $R _{1} ^{\prime
}/R ^{\prime }$;
\par
{\rm (c)} There exists is a sequence $R _{n}$, $n \in \mathbb{N}$, of
finite extensions of $K$ in $R$, such that $p ^{n}$ divides the
degree $[R _{n}\colon K]$, for each index $n$.
\end{lemm}
\par
\medskip
\begin{proof}
Our starting point is the fact that $K$ is a quasilocal field, in the
sense of \cite{Ch1}, which implies the following:
\par
\medskip\noindent
(3.6) For any pair $\Lambda , \Lambda ^{\prime }$ of finite
extensions of $K$, such that $\Lambda \in I(\Lambda ^{\prime }/K)$
and $\Lambda ^{\prime }$ is separable over $\Lambda $, the group
Br$(\Lambda ^{\prime }/\Lambda )$ consists of all elements of
Br$(\Lambda )$ of orders dividing $[\Lambda ^{\prime }\colon \Lambda
]$ (see \cite{S}, Ch. XIII, Sect. 3; \cite{Ch1}, Corollary~8.5); the
assertion holds without the assumption that $\Lambda ^{\prime
}/\Lambda $ is separable, provided that $K$ is an almost perfect
field (cf. \cite{Ch1}, Corollary~8.6).
\par
\medskip\noindent
Note also that Br$(R) _{p}$ equals the union of the images of
Br$(\Lambda ) _{p}$ under the scalar extension maps $\pi _{\Lambda
\to R}$, when $\Lambda $ runs across the set $I _{0}(R/K)$ of finite
extensions of $K$ in $R$ (cf., e.g., \cite{Ch1}, (1.3)). This
indicates that Br$(R) _{p} = \{0\}$ if and only if Br$(\Lambda ) _{p}
\subseteq {\rm Br}(R/\Lambda )$, for any $\Lambda \in I _{0}(R/K)$.
Observe now that if $R _{n}$, $n \in \mathbb{N}$, are fields
satisfying condition (c) of Lemma \ref{lemm3.4} with respect to
$R/K$, then for each finite extension $\Lambda $ of $K$ in $R$, the
sequence $\Lambda R _{n}$, $n \in \mathbb{N}$, has an infinite
subsequence satisfying the same condition with respect to $R/\Lambda
$. At the same time, the violation of condition (c) means that there
is a finite extension $\Lambda _{1}$ of $K$ in $R$, such that $p$
does not divide the degree of any finite extension of $\Lambda _{1}$
in $R$. Therefore, by (3.6), $\pi _{\Lambda _{1} \to R}$ maps
Br$(\Lambda _{1}) _{p}$ injectively into Br$(R) _{p}$. Since Br$(K)$ is
isomorphic to the quotient group $\mathbb{Q}/\mathbb{Z}$ (cf.
\cite{S}, Ch. XIII, Sect. 3), whence, Br$(L) _{p} \neq \{0\}$ for
every finite extension $L/K$, these remarks enable one to deduce from
(3.6) (under the hypothesis that $R _{n}/K$, $n \in \mathbb{N}$, are
separable, or $K$ is almost perfect) the following: (i) conditions (c)
and (a) of Lemma \ref{lemm3.4} are equivalent; (ii) condition (c)
holds if and only if Br$(R) _{p} = \{0\}$. For a proof of the
equivalence of conditions (a) and (b), and of the former part of
condition (a) to the latter one, we refer the reader to \cite{S1},
Ch. II, 3.1, and \cite{GiSz}, Theorem~6.1.8.
\par
It remains to be seen that if $p \neq {\rm char}(K)$, then the
assertions of the lemma hold, for any algebraic extension $R/K$.
Let $R _{0}$ be the maximal separable extension of $K$ in $R$. It
is well-known that if $R \neq R _{0}$, then char$(K) = q > 0$ and
finite extensions of $R _{0}$ in $R$ are of $q$-primary degrees
(cf. \cite{L}, Ch. V, Corollary~6.2). This enables one to prove that
Br$(R/R _{0})$ is a subgroup of Br$(R _{0}) _{q}$ (see \cite{P},
Sect. 13.4). Note finally that, by the Albert-Hochschild theorem (cf.
\cite{S1}, Ch. II, 2.2), $\pi _{R _{0} \to R}$ maps Br$(R _{0}) _{p}$
surjectively upon Br$(R) _{p}$. When $p \neq q$, these
observations show that $\pi _{R _{0} \to R}$ induces an
isomorphism Br$(R _{0}) _{p} \cong {\rm Br}(R) _{p}$, so we have Br$(R
_{0}) _{p} = \{0\} \Leftrightarrow {\rm Br}(R) _{p} = \{0\}$. The
obtained equivalence completes the proof of Lemma
\ref{lemm3.4}, since its assertions apply to the extension $R
_{0}/K$.
\end{proof}
\par
\medskip
It is known that, for any field $R$, conditions (a) and (b) stated
in Lemma \ref{lemm3.4} are equivalent, and when they hold,
$\mathcal{G}_{R}$ is a profinite group of cohomological
$p$-dimension cd$_{p}(\mathcal{G}_{R}) \le 1$; this implication is
an equivalence in case $R$ is perfect or $p \neq {\rm char}(R)$
(cf. \cite{S1}, Ch. II, 3.1, and \cite{GiSz}, Theorem~6.1.8). Thus
it follows that dim$(R) \le 1 \Leftrightarrow N(R _{1} ^{\prime
}/R ^{\prime }) = R ^{\prime \ast }$, for any pair $R _{1}
^{\prime } \in {\rm Fe}(R)$, $R ^{\prime } \in I(R _{1} ^{\prime
}/R)$ (cf. \cite{S1}, Ch. II, 3.1). Clearly, $N(R _{1} ^{\prime
}/R ^{\prime }) = R ^{\prime \ast }$ if and only if the $R
^{\prime }$-form $N _{R _{1}'/R'}(X _{1}, \dots , X _{n}) - aX
_{n+1} ^{n}$ has a nontrivial zero over $R ^{\prime }$, for each
$a \in R ^{\prime \ast }$, where $n = [R _{1} ^{\prime }\colon R
^{\prime }]$ and $N _{R _{1}'/R'}(X _{1}, \dots , X _{n})$ is the
norm form of degree $n$ in algebraically independent variables $X
_{1}, \dots , X _{n}$ over $R ^{\prime }$, associated with a fixed
$R ^{\prime }$-basis of $R _{1} ^{\prime }$. This standardly
proves the fact that $C _{1}$-fields have dimension $\le 1$ (cf.
\cite{S1}, Ch. II, 3.2) and leads to the problem of whether the
converse holds if one restricts to fields from special classes of
sufficient research interest, such as the class of Henselian
fields algebraic over a global field $K$, and the class
$\mathcal{V}(K)$ defined in Section 2 (for other classes, see
\cite{Koe}, Theorem~B, and \cite{N}).
\par
\medskip
\begin{lemm}
\label{lemm3.5}
Let $(K, v)$ be an {\rm HDV}-field with $\widehat K$ quasifinite, and
let $R/K$ be an algebraic extension. Then {\rm Br}$(R) _{p} \neq
\{0\}$, for a given $p \in \mathbb{P}$, if and only if $v(R) \neq
pv(R)$ and $\widehat R(p) \neq \widehat R$.
\end{lemm}
\par
\medskip
\begin{proof}
Suppose first that $\widehat R(p) \neq \widehat R$ and $v(R) \neq
pv(R)$. Then, by Lemma \ref{lemm3.1}, there exists a degree $p$
extension $R _{1}$ of $R$ in $R _{\rm ur}$; also, by assumption,
there is $\pi \in K ^{\ast }$ of value $v(\pi ) \notin pv(R)$. Fix a
generator $\sigma $ of $\mathcal{G}(R _{1}/R)$ and consider the
cyclic $R$-algebra $A = (R _{1}/R, \sigma , \pi )$ (for the
definition of $A$, see, e.g., \cite{GiSz}, Sect. 2.5). It is known
that $A$ is a nicely semiramified division $R$-algebra of dimension
$p ^{2}$, in the sense of Jacob-Wadsworth (see \cite{TW}, page 452,
and further references there). In particular, $R$ equals the centre
of $A$ and the Brauer equivalence class of $A$ is an element of order
$p$ in Br$(E)$, which shows that Br$(R) _{p} \neq \{0\}$.
\par
In the rest of the proof of Lemma \ref{lemm3.5}, we assume that
Br$(R) _{p} \neq \{0\}$. Our goal is to show that $\widehat R(p)
\neq \widehat R$ and $v(R) \neq pv(R)$. We first prove that $\widehat
R(p) \neq \widehat R$. As $\widehat K$ is quasifinite,
$\mathcal{G}_{\widehat K}$ is isomorphic to the profinite group
$\prod _{\ell \in \mathbb{P}} \mathbb{Z} _{\ell }$, so it follows
from Lemma \ref{lemm3.1} and the equality $R _{\rm ur} = RK _{\rm
ur}$, that if $\widehat R(p) = \widehat R$, then $R \cap K _{\rm ur}$
contains as a subfield a $\mathbb{Z} _{p}$-extension $\Gamma _{p}$ of
$K$. Hence, by Galois theory, condition (c) of Lemma \ref{lemm3.4}
holds with respect to $\Gamma _{p}/K$. This requires Br$(\Gamma _{p}) _{p} =
\{0\}$ and Br$(R) _{p} = \{0\}$, which is a contradiction proving
that $\widehat R(p) \neq \widehat R$.
\par
We turn to the proof of the assertion that $v(R) \neq pv(R)$. Denote
by $R _{0}$ the maximal separable extension of $K$ in $R$. Since
Br$(R) _{p} \neq \{0\}$, Lemma \ref{lemm3.4} yields Br$(R _{0}) _{p}
\neq \{0\}$ and implies the existence of a finite extension $L$ of
$K$ in $R _{0}$, such that $p \nmid [L _{1}\colon L]$, for any finite
extension $L _{1}$ of $L$ in $R _{0}$. In view of (3.3) (or the
isomorphism $v(L) \cong \mathbb{Z}$, see \cite{E3},
Corollary~14.2.2), this means that $v(R _{0}) \neq pv(R _{0})$, which
completes the proof of Lemma \ref{lemm3.5} in the case where char$(K)
= 0$. We assume further that char$(K) = q > 0$. Then $R$ is purely
inseparable over $R _{0}$, whence finite extensions of $R _{0}$ in
$R$ are of $q$-primary degrees. Therefore, (3.3) indicates that if $p
\neq q$, then $v(R) \neq pv(R)$, as claimed.
\par
Suppose now that $p = q$, put $R _{v} = R _{v_{R}}$ and $L _{v} = L
_{v _{L}}$, fix an algebraic closure $\overline R _{v}$ of $R _{v}$,
denote by $\overline R$ the algebraic closure of $R$ in $\overline R
_{v}$, identify the field $L _{v}$ with the topological closure of
$L$ in $R _{v}$, and put $v' = \overline {v _{L}}$. Consider the fields $R
^{\prime } = RL _{v}$ and $R _{0} ^{\prime } = R _{0}L _{v}$. We
first show that $R ^{\prime } = R _{0} ^{\prime }L ^{\prime }$, where
$L ^{\prime }$ is the maximal purely inseparable extension of $L
_{v}$ in $R ^{\prime }$. Our proof  relies on the fact that, by (3.2)
(b), $(L _{v}, v')$ is a complete discrete valued field with
$\widehat L _{v} = \widehat L$. Observing also that $\widehat
L/\widehat K$ is a finite extension, whence $\widehat L$ is
a quasifinite field, one obtains from (3.5) (a) (or rather, from the
fulfillment of condition (3.4) (b)) that $[L _{v}\colon L _{v} ^{p}]
= p$. As the extension $R _{0} ^{\prime }/L _{v}$ is separable (it
preserves the separability of $R _{0}/L$), this allows to deduce from
\cite{BH}, Lemma~2.12 (and \cite{L}, Ch. V, Corollary~6.10) that $[R
_{0} ^{\prime }\colon R _{0} ^{\prime p}] = p$. Thus it turns out
that, for each $\nu \in \mathbb{N}$, there exists a unique pair
$\Lambda _{\nu }, \Lambda _{\nu } ^{\prime } \in I(\overline R _{v}/L
_{v})$, such that $\Lambda _{\nu }/L _{v}$ and $\Lambda _{\nu }
^{\prime }/R _{0} ^{\prime }$ are purely inseparable extensions of
degree $p ^{\nu }$. Moreover, one obtains the following:
\par
\medskip\noindent
(3.7) (a) $\Lambda _{1} ^{p} = L _{v}$, and for each $\nu \in
\mathbb{N}$, $\Lambda _{\nu } \subset \Lambda _{\nu +1}$ and $\Lambda
_{\nu +1} ^{p} = \Lambda _{\nu }$; also, for any $\pi \in L _{v}
^{\ast p}$ with $v'(\pi ) \notin pv'(L _{v}) = pv(L)$, $\Lambda _{\nu
} = L _{v}(\pi _{\nu })$, where $\pi _{\nu } \in \overline R _{v}$ is
the $p ^{\nu }$-th root of $\pi $;
\par
(b) The union $\Lambda _{\infty }$ of the fields $\Lambda _{\nu
}$, $\nu \in \mathbb{N}$, is a perfect field, and $I(\Lambda _{\infty
}/L _{v})$ equals the set $\{\Lambda _{\infty }, L _{v}, \Lambda
_{\nu }\colon \nu \in \mathbb{N}\}$ (this description of $I(\Lambda
_{\infty }/L _{v})$ characterizes the property that $[L _{v}\colon L
_{v} ^{p}] = p$, see Remark~1 in \cite{S1}, Ch. II, 3.2).
\par
\medskip\noindent
Clearly, $\Lambda _{\infty }$ is the maximal purely inseparable
extension of $L _{v}$ in $\overline R _{v}$. As $R _{0} ^{\prime }$
is separable over $L _{v}$, this implies the set of purely
inseparable extensions of $R _{0} ^{\prime }$ in $\overline R _{v}$
is equal to $I(\Lambda _{\infty } ^{\prime }/R _{0} ^{\prime })$,
where $\Lambda _{\infty } ^{\prime } = R _{0} ^{\prime }\Lambda
_{\infty }$ (cf. \cite{L}, Ch. V, Corollary~6.10), and for each $\nu
\in \mathbb{N}$, $R _{0} ^{\prime }\Lambda _{\nu }/R _{0} ^{\prime }$
is a field extension of degree $p ^{\nu }$. It is now easy to see
that $\Lambda _{\nu } ^{\prime } = R _{0} ^{\prime }\Lambda _{\nu }$,
$\nu \in \mathbb{N}$, and also, $R _{0} ^{\prime }$ and
$\Lambda _{\nu } ^{\prime }$, $\nu \in \mathbb{N} \cup \{\infty \}$,
are all purely inseparable extensions of $R _{0} ^{\prime }$ in
$\overline R _{v}$, so the equality $R ^{\prime } = R _{0} ^{\prime
}L ^{\prime }$ is proved.
\par
We continue with the proof of the inequality $v(R) \neq pv(R)$. As
noted above, the field extension $R _{0} ^{\prime }/L _{v}$ is
separable, and since $R ^{\prime } = R _{0} ^{\prime }L ^{\prime
}$, it follows that $R ^{\prime }/L ^{\prime }$ is also separable
whereas $R ^{\prime }/R _{0} ^{\prime }$ is purely inseparable.
Identifying the completion of $(R _{0}, v _{R_{0}})$ with the
topological closure $R _{0,v}$ of $R _{0}$ in $R _{v}$, one
obtains similarly that $R _{0,v}.R/R _{0,v}$ is purely inseparable
as well. Observe that $(L, v _{L})$, $(L _{v}, v')$ and $R _{0}/L$
satisfy the conditions of Lemma \ref{lemm3.3}. Note also $(R _{v},
\overline {v _{R}})/(L _{v}, v')$ is a valued field extension, so
it follows from Lemma \ref{lemm3.3}, applied to $(L, v _{L})$, $(L
_{v}, v')$ and $R _{0}$, and from the uniqueness of the
prolongations of $v _{R _{0}}$, $v'_{R _{0}'}$ and $\overline {v
_{R _{0}}}$ on
 $R$, $R ^{\prime }$, and $R _{0,v}R$, respectively (cf. \cite{E3},
Proposition~14.2.5), that $(R ^{\prime }, v'_{R'})$ is an
intermediate valued field of $(R _{v}, \overline {v _{R}})/(R, v
_{R})$. This implies the topologies on $R$ and $R ^{\prime }$
associated with $v _{R}$ and $v' _{R'}$, respectively, are induced
by the topology of $R _{v}$ (which in turn is determined by
$\overline {v _{R}}$, see \cite{E3}, Theorem~9.3.2). Moreover, $R$
and $R ^{\prime }$ are dense in $R _{v}$, and both $(R, v _{R})$
and $(R ^{\prime }, v'_{R'})$ are Henselian, which guarantees the
injectivity of the maps $\pi _{R \to R _{v}}$ and $\pi _{R' \to R
_{v}}$ (cf. \cite{Cohn}, Theorem~1). Since, by basic well-known
properties of tensor products (cf. \cite{P}, Sect. 9.4,
Corollary~a), $\pi _{R \to R _{v}}$ equals the composition $\pi
_{R' \to R _{v}} \circ \pi _{R\to R'}$, this implies $\pi _{R \to
R'}$ is injective, so it follows from the nontriviality of Br$(R)
_{p}$ that Br$(R ^{\prime }) _{p} \neq \{0\}$. Note further that
$[L ^{\prime }\colon L _{v}] < \infty $. Assuming the opposite,
one obtains from (3.7) (b) that $L ^{\prime }$ must equal the
(perfect) field $\Lambda _{\infty }$. Therefore, $R ^{\prime } = R
_{0} ^{\prime }L ^{\prime }$ must also be perfect, which requires
Br$(R ^{\prime }) _{p} = \{0\}$ (cf. \cite{A1}, Ch. VII,
Theorem~22) - a contradiction proving that $[L ^{\prime }\colon L
_{v}] < \infty $.
\par
It is now easy to complete the proof of the inequality $v(R) \neq
pv(R)$ (and of Lemma \ref{lemm3.5}). Indeed, $(L _{v}, v')$ is an
HDV-field, and we have $[L ^{\prime }\colon L _{v}] < \infty $, so
it follows from \cite{E3}, Corollary~14.2.2, that $(L ^{\prime },
v'_{L'})$ is HDV as well. As $[L _{v}\colon L _{v} ^{p}] = p$ and
$L ^{\prime }/L$ is purely inseparable, (3.7) (a) shows that $L
^{\prime }/L$ is totally ramified; in particular, $\widehat L
^{\prime } = \widehat L _{v} = \widehat L$. Note further that
$\widehat L$ is a quasifinite field, $R ^{\prime }/L ^{\prime }$
is a separable extension and Br$(R ^{\prime }) _{p} \neq \{0\}$.
Therefore, $R ^{\prime }/L ^{\prime }$ violates condition (c) of
Lemma \ref{lemm3.4}, which enables one to deduce from (3.3) and
the inequality $v'(L ^{\prime }) \neq pv'(L ^{\prime })$ that
$v'(R ^{\prime }) \neq pv'(R ^{\prime })$. Since $(R ^{\prime },
v' _{R'})$ is an intermediate valued field of $(R _{v}, \overline
{v _{R}})/(R, v_{R})$, whence, $\overline {v _{R}}(R _{v}) = v'(R
^{\prime }) = v(R)$, the inequality $v(R) \neq pv(R)$ now becomes
obvious, so Lemma \ref{lemm3.5} is proved.
\end{proof}
\par
\medskip
Lemma \ref{lemm3.5} shows that if $(K, v)$ is an HDV-field with
$\widehat K$ quasifinite, and $R/K$ is an algebraic extension, then
dim$(R) \le 1$ if and only if the intersection $S(R) \cap \Sigma (R)$
is the empty set, where $S(R) = \{p \in \mathbb{P}\colon
\mathcal{G}(\widehat R(p)/\widehat R) \cong \mathbb{Z} _{p}\}$ and
$\Sigma (R) = \{p \in \mathbb{P}\colon v(R) \neq pv(R)\}$.
\par
\medskip
\begin{rema}
\label{rema3.6} Let $(K, v)$ be an {\rm HDV}-field with $\widehat K$
quasifinite and {\rm char}$(K) = p > 0$, and let $R/K$ be any
algebraic field extension. Then it follows from (3.3) and Lemma
\ref{lemm3.5} that {\rm Br}$(R) _{p} = \{0\}$ if and only if one of
the conditions {\rm (a)} and {\rm (b)} of Lemma \ref{lemm3.4} holds;
also, by the proof of Lemma \ref{lemm3.4}, the assumption that {\rm
Br}$(R) _{p} = \{0\}$ ensures that condition {\rm (c)} is satisfied.
The fulfillment of condition {\rm (c)}, however, does not guarantee
that {\rm Br}$(R) _{p} = \{0\}$ if $K$ has an infinite purely
inseparable extension $K ^{\prime }$ with $v(K ^{\prime }) = v(K)$
and $\widehat K ^{\prime } = \widehat K$ (see \cite{Ch1}, Remark~8.7,
for examples of such HDV-fields). Then $(K ^{\prime }, v _{K'})$ is
{\rm HDV} and $\widehat K ^{\prime }$ is quasifinite, proving that
{\rm Br}$(K ^{\prime }) \cong \mathbb{Q}/\mathbb{Z}$; hence, {\rm
Br}$(K ^{\prime }) _{p'} \neq \{0\}$, $p' \in \mathbb{P}$. On the
other hand, there are fields $R _{n} \in I(K ^{\prime }/K)$, $n \in
\mathbb{N}$, such that $[R _{n}\colon K] = p ^{n}$, for each $n$;
therefore, $K ^{\prime }/K$ satisfies condition {\rm (c)} of Lemma
\ref{lemm3.4}.
\end{rema}
\par
\medskip
At the end of this Section, note that, for an arbitrary field $R$, the
assumption that Br$(R) _{p} = \{0\}$ is, generally, weaker than
conditions (a) and (b) of Lemma \ref{lemm3.4}. As shown by M.
Auslander, the formal power series field $\mathbb{Q} _{\rm
sol}((X))$, where $\mathbb{Q} _{\rm sol}$ is the compositum of finite
Galois extensions of $\mathbb{Q}$ in $\mathbb{Q} _{\rm sep}$ with
solvable Galois groups, satisfies Br$(\mathbb{Q} _{\rm sol}((X))) =
\{0\}$ but violates condition (a), for each $p \in \mathbb{P}$ (see
\cite{S1}, Ch. II, 3.1). If, however, $E \in I(\mathbb{Q} _{\rm
sep}/\mathbb{Q})$ and Br$(E) _{p} = \{0\}$, for some $p \in
\mathbb{P}$, then Br$(E ^{\prime }) _{p} = \{0\}$ whenever $E
^{\prime } \in I(\mathbb{Q} _{\rm sep}/E)$ (see \cite{FS},
Theorem~4); this can also be deduced from Lemma \ref{lemm3.4} and
\cite{S1}, Ch. II, Proposition~9.
\par
\medskip
\section{\bf Proof of Theorem \ref{theo2.2}}
\par
\medskip
Let $(K, v)$ be an HDV-field with $\widehat K$ quasifinite, and let
$L/K$ be an algebraic extension, $S(L) = \{p \in \mathbb{P}\colon
\mathcal{G}(\widehat L(p)/\widehat L) \cong \mathbb{Z} _{p}\}$ and
$\Sigma (L) = \{p \in \mathbb{P}\colon v(L) \neq pv(L)\}$. Then, by
Lemma \ref{lemm3.5}, dim$(L) \le 1$ if and only if $S(L) \cap \Sigma
(L) = \emptyset $. Moreover, by Lang's theorem (see \cite{L1},
Theorem~10), $L$ is a $C _{1}$-field if $S(L) = \emptyset $ and
$K _{v} = K$ (the emptiness of $S(L)$ ensures that $K _{\rm ur} \in
I(L/K)$, whence, $L$ preserves the $C _{1}$ type of $K _{\rm ur}$).
Therefore, in this Section, we prove Theorem \ref{theo2.2}
considering fields $L$ with dim$(L) \le 1$ and $S(L) \neq \emptyset
$. Our proof relies on the following lemma:
\par
\medskip
\begin{lemm}
\label{lemm4.1}
Let $(K, v)$ be an {\rm HDV}-field with $\widehat K$ quasifinite, and
$S$, $\Sigma $ be nonempty proper subsets of
$\mathbb{P}$. Then there exists an algebraic extension $E/K$, such
that $S(E) = S$ and $\Sigma (E) = \Sigma $; hence, {\rm Br}$(E) _{p}
\neq \{0\}$ if and only if $p \in S \cap \Sigma $.
\end{lemm}
\par
\medskip
\begin{proof}
Denote by $\widehat K(S)$ the compositum of the maximal
$p$-extensions $\widehat K(p)$, $p \in \mathbb{P} \setminus S$, of
$\widehat K$ in $\widehat K _{\rm sep}$. It follows from Lemma
\ref{lemm3.1} that $K$ has a Galois extension $K(S)$ in $K _{\rm ur}$
with $\widehat {K(S)} = \widehat K(S)$ and $\mathcal{G}(K(S)/K) \cong
\mathcal{G}(\widehat K(S)/\widehat K)$. As $\widehat K$ is
quasifinite and $K(S) _{\rm ur} = K _{\rm ur}$, this enables one to
deduce from Galois theory (cf. \cite{L}, Ch. VI, Theorem~1.12) and
Lemma \ref{lemm3.1} that $\mathcal{G}_{\widehat {K(S)}}$ is
isomorphic to the topological group product $\prod _{p \in S}
\mathbb{Z} _{p}$. Moreover, Br$(K(S))_{p} = \{0\}$, $p \in \mathbb{P}
\setminus S$, by Lemma \ref{lemm3.5}. Now fix an algebraic closure
$\overline K$ of $K _{\rm sep}$, take a generator $\pi $ of the ideal
$M _{v}(K)$, and for each $t \in \mathbb{P} \setminus \Sigma $,
denote by $\Theta _{t}$ the set $\{\pi _{t,n}\colon n \in
\mathbb{N}\}$, where $\pi _{t,n} \in \overline K\colon n \in
\mathbb{N}$, is a sequence defined inductively so that $\pi _{t,1} =
\pi $, and $\pi _{t,n} ^{t} = \pi _{t,(n-1)}$, for every $n \ge 2$.
It is easily verified that for any $t \in \mathbb{P} \setminus \Sigma
$, the extension $K(\Theta _{t})/K$ is infinite and finite
subextensions of $K$ in $K(\Theta _{t})$ are totally ramified of
$t$-primary degrees. Therefore, by Lemma \ref{lemm3.5}, Br$(K(\Theta
_{t}))_{t} = \{0\}$ and $v(K(\Theta _{t})) = tv(K(\Theta _{t}))$, for
$t \in \mathbb{P} \setminus \Sigma $.
\par
Consider now the compositum $E$ of the fields $K(S)$ and $K(\Theta
_{t})$, $t \in \mathbb{P} \setminus \Sigma $. The described
properties of $K(S)$ and $K(\Theta _{t})$, $t \in \mathbb{P}
\setminus \Sigma $, ensure that $(E, v)$ satisfies the following
conditions:
\par
\medskip\noindent
(4.1) (a) Finite extensions of $K(S)$ in $E$ are totally ramified
and their degrees are not divisible by any $\lambda \in \Sigma $;
\par
(b) $\mathcal{G}_{\widehat E}$ is isomorphic to $\prod _{p \in S}
\mathbb{Z} _{p}$, and $v(E) \neq pv(E)$ if and only if $p \in \Sigma
$.
\par
\medskip\noindent
Statements (4.1) imply $S(E) = S$ and $\Sigma (E) = \Sigma $, so
it follows from Lemma \ref{lemm3.5} that Br$(E) _{p} \neq \{0\}$ if
and only if $p \in S \cap \Sigma $. Lemma \ref{lemm4.1} is proved.
\end{proof}
\par
\medskip
Next we show that $\mathbb{P}$ possesses subsets $S$ and $\Sigma $ of
satisfying the following:
\par
\medskip\noindent
(4.2) (a) $\Sigma $ is infinite, $S \cap \Sigma = \emptyset $, and
for each $p \in \Sigma $, there is $k(p) \in \mathbb{Z}$, such that
$2 \le k(p) \le (p - 1)/2$ and all prime divisors of $k(p).(p -
k(p))$ lie in $S$;
\par
(b) $\gcd (k(p).(p-k(p)), k(p').(p'-k(p')))$ is a $2$-primary
number, for any
\par\noindent
$p, p' \in \Sigma $ with $p \neq p'$.
\par
\medskip\noindent
The proof of (4.2) relies on Dirichlet's theorem about the prime
numbers in an arithmetic progression, and on the Chinese Remainder
Theorem. Using repeatedly these theorems, one obtains inductively
that there exist positive integers $k _{n}, p _{n}$, $n \in \mathbb
N$, such that:
\par
\medskip\noindent
(4.3) (a) $k _{1} = 2$, $p _{1} = 5$, and for each $n$, $k _{n}
\equiv 2 ^{n} ({\rm mod} \ 2 ^{n+1})$ and $p _{n} \in \mathbb{P}$;
\par
(b) If $n \ge 2$, then $k _{n} \equiv 1 \ ({\rm mod} \ \prod _{j=1}
^{n-1} p _{j}k _{j}.2 ^{-j}(p _{j} - k _{j}))$,
\par
\noindent
$p _{n} \equiv 2 \ ({\rm mod}\ \prod _{j=1} ^{n-1} p _{j}k _{j}.2
^{-j}(p _{j} - k _{j}))$, $p _{n} \equiv 1 ({\rm mod}\ k _{n})$, and
$p _{n} \ge 1 + 2k _{n}$.
\par
\medskip\noindent
Let now $\Sigma = \{p _{n}\colon n \in \mathbb{N}\}$ and $S = \{p
\in \mathbb{P}\colon p \mid k _{n}(p _{n} - k _{n}),$ for some $n
\in \mathbb{N}\}$. Arguing by induction on $n$, and using (4.3), one
obtains that $p _{n} \notin S$, for any $n$, i.e. $\Sigma \cap S =
\emptyset $. Moreover, it is easily verified that $S$ and $\Sigma $
satisfy conditions (4.2), where $k(p _{n}) = k _{n}$, for each $n
\in \mathbb{N}$.
\par
\medskip
We prove Theorem \ref{theo2.2}. Suppose that $S$ and $\Sigma $ are
defined in accordance with (4.3), and $(K, v)$ and $E$ have the
properties required by Lemma \ref{lemm4.1} (equivalently, by (4.1)
(b)). Fix elements $\pi _{n} \in E ^{\ast }$, $n \in \mathbb{N}$, of
values $v(\pi _{n}) > 0$ and $v(\pi _{n}) \notin p _{n}v(E)$, for any
$n$. It follows from (4.1) (b), Galois theory (cf. \cite{L}, Ch. VI)
and the definition of $S$ that $E$ has cyclic extensions $E _{n}$ and
$E _{n} ^{\prime }$, such that $[E _{n}\colon E] = k _{n}$, $[E _{n}
^{\prime }\colon E] = p _{n} - k _{n}$ and $E _{n}.E _{n} ^{\prime }
\subset E _{\rm ur}$, for each $n \in \mathbb{N}$. Take primitive
elements $\xi _{n} \in O _{v}(E _{n})$ of $E _{n}/E$ and $\eta _{n}
\in O _{v}(E _{n} ^{\prime })$ of $E _{n} ^{\prime }/E$, so that the
residue classes $\hat \xi _{n} \in \widehat E _{n}$ and $\hat \eta
_{n} \in \widehat E _{n} ^{\prime }$ are primitive elements of
$\widehat E _{n} ^{\prime }/\widehat E$ and $\widehat E _{n} ^{\prime
}/\widehat E$, respectively, and consider a system $\overline X _{n}
= X _{n,i}$, $i = 1, \dots , k _{n}$, of algebraically independent
variables over $E$. In view of Lemma \ref{lemm3.1} and Galois theory,
$E _{n}(\overline X _{n})/E(\overline X _{n})$ and $E _{n} ^{\prime
}(\overline X _{n})/E(\overline X _{n})$ are cyclic field extensions
of degrees $k _{n}$ and $p _{n} - k _{n}$, respectively, and the
product $g_{n}(\overline X _{n})$ of the norms
$$N _{E _{n}(\overline X _{n})/E(\overline X _{n})}(\sum _{i=1}
^{k_{n}} \xi _{n} ^{i-1}X _{n,i}) \ {\rm and} \ N _{E _{n}'(\overline
X _{n})/E(\overline X _{n})}(\sum _{i=1} ^{k _{n}} \eta _{n} ^{i-1}X
_{n,i})$$
is a form of degree $p _{n}$ in the variables $X _{n,i}$, $i = 1,
\dots , k _{n}$, with coefficients in $O _{v}(E)$. Moreover, since
$v$ is Henselian, it follows from the inclusion $E_{n}E _{n} ^{\prime
} \subset E _{\rm ur}$, the assumptions on $\xi _{n}$ and $\eta
_{n}$, and the inequality $k _{n} <p _{n} - k _{n}$ that
\par
\medskip\noindent
(4.4) $g _{n}(\bar \alpha _{n}) \in O _{v}(E)$, for every $k
_{n}$-tuple $\bar \alpha _{n} = (\alpha _{n,1}, \dots , \alpha
_{n,k_{n}})$ with components $\alpha _{n,i} \in O _{v}(E)$; $g
_{n}(\bar \alpha _{n}) \in M _{v}(E)$ if and only if $\alpha _{n,i}
\in M _{v}(E)$, $i = 1, \dots , k _{n}$.
\par
\medskip\noindent
One also deduces from (4.4) that $v(g _{n}(\bar \beta _{n})) \in p
_{n}v(E)$, provided that
\par\noindent
$\bar \beta _{n} = (\beta _{n,1}, \dots , \beta _{n,k_{n}})$,
$\beta _{n,i} \in E$, for $i = 1, \dots k _{n}$, and $\beta _{n,i'}
\neq 0$, for some index $i'$. Let now $\widetilde X = X
_{n,i,j}\colon $ $i = 1, \dots , k _{n}$; $j = 1, \dots , p _{n}$, be
a system of $p _{n}k _{n}$ algebraically independent variables over
$E$, and let $f _{n}(\widetilde X) = \sum _{j=1} ^{p _{n}} g
_{n}(\widetilde X _{n,j}).\pi _{n} ^{j}$, where $\widetilde X _{n,j}
= X _{n,i,j}$, $i = 1, \dots , k _{n}$, for each fixed $j \in \{1,
\dots , p _{n}\}$. Using the noted properties of $g _{n}(\overline
X)$ and the fact that $v(\pi _{n}) \notin p _{n}v(E)$, one obtains
that
\par
\medskip\noindent
(4.5) For any $n \in \mathbb{N}$, $f _{n}(\widetilde X)$ is an
$E$-form of degree $p _{n}$ in $p _{n}k _{n}$ variables, and without
a nontrivial zero over $E$; in particular, $E$ is not a $C
_{1}$-field.
\par
\medskip\noindent
It follows from (4.1) (b) and Galois theory that if $E ^{\prime }$ is
a finite extension of $E$, then $S(E ^{\prime }) = S(E)$ and $\Sigma
(E ^{\prime }) = \Sigma (E)$. Therefore, $E ^{\prime }$ preserves the
properties of $E$ described by (4.1) (b), which enables one to show
by the method of proving (4.5) that $E ^{\prime }$ is not a $C
_{1}$-field. Theorem \ref{theo2.2} is proved.
\par
\medskip
Theorem \ref{theo2.2} and our next result exhibit the fact that, for
any global or local field $K$, algebraic extensions $E$ of $K$ with
dim$(E) \le 1$ need not be $C _{1}$-fields.
\par
\medskip
\begin{coro}
\label{coro4.2}
Let $K$ be a global field, $v$ a Krull valuation of $K$, and $K
^{\prime }$ the maximal separable extension of $K$ in $K _{v}$. Then
$K ^{\prime }$ has an algebraic extension $E ^{\prime }$, such that
{\rm dim}$(E ^{\prime }) \le 1$ and finite extensions of $E ^{\prime
}$ are not $C _{1}$-fields.
\end{coro}
\par
\medskip
\begin{proof}
Denote by $\bar v$ the continuous prolongation of $v$ on $K _{v}$,
and by $v ^{\prime }$ the valuation of $K ^{\prime }$ induced by
$\bar v$. Observe that $(K ^{\prime }, v ^{\prime })$ is HDV as well
as an intermediate valued field of $(K _{v}, \bar v)/(K, v)$ (see (3.2)
and the end of Section~1). Hence, $v(K ^{\prime }) = v(K)$ and
$\widehat K ^{\prime }$ equals the residue field of $(K, v)$; in
addition, $\widehat K ^{\prime }$ is finite. Now Corollary
\ref{coro4.2} is proved by applying Theorem \ref{theo2.2} to $(K
^{\prime }, v ^{\prime })$.
\end{proof}
\par
\medskip
\begin{coro}
\label{coro4.3}
Let $(K, v)$ be an {\rm HDV}-field with $\widehat K$ quasifinite,
and let $E$ be an algebraic extension of $K$ with $S(E) \cap \Sigma
(E) = \emptyset $. Assume that $k _{1}$ and $k _{2}$ are integers,
such that $2 \le k _{1} < k _{2}$, every divisor $p \in \mathbb{P}$
of $k _{1}k _{2}$ lies in $S(E)$, and every $p' \in \mathbb{P}$
dividing $k _{1} + k _{2}$ lies in $\Sigma (E)$. Then {\rm dim}$(E)
\le 1$ and there exists an $E$-form $f$ of degree $k _{1} + k _{2}$
in $k _{1}(k _{1} + k _{2})$ variables, without a nontrivial zero
over $E$; in particular, $E$ is not a $C _{1}$-field and $S(E)$
contains at least $2$ elements.
\end{coro}
\par
\medskip
\begin{proof}
The inequality dim$(E) \le 1$ follows from Lemma \ref{lemm3.5}
and the condition on $S(E) \cap \Sigma (E)$. Also, our assumptions
show that $\gcd (k _{1}, k _{2}) = 1$, $S(E)$ is of cardinality $\ge
2$, and there exist fields $E _{j} \in I(E _{\rm ur}/E)$ with $[E
_{j}\colon E] = k _{j}$, for $j = 1, 2$. Using the cyclicity of
finitely-generated subgroups of $\mathbb{Q}$, one obtains that the
set $(k _{1} + k _{2})v(E) = \{(k _{1} + k _{2})\gamma \colon \gamma
\in v(E)\}$ is a subgroup of $v(E)$ of index $k _{1} + k _{2}$, and
the group $v(E)/(k _{1} + k _{2})v(E)$ is cyclic. This allows to
define (like in the proof of (4.5)) an $E$-form as required by
Corollary \ref{coro4.3}.
\end{proof}
\par
\medskip
Note that if $k _{1}$ and $k _{2}$ are integers with $2 \le k _{1} <
k _{2}$ and $\gcd (k _{1}, k _{2}) = 1$, then $\gcd (k _{1}k _{2}, k
_{1} + k _{2}) = 1$. This implies the existence of various subsets
$S$ and $\Sigma $ of $\mathbb{P}$, such that $S \cap \Sigma =
\emptyset $, and $p \in S$, $p' \in \Sigma $, for any pair $p, p' \in
\mathbb{P}$ satisfying $p \mid k _{1}k _{2}$ and $p' \mid k _{1} + k
_{2}$. Therefore, Lemma \ref{lemm4.1}, Corollary \ref{coro4.3} and
the proof of Theorem \ref{theo2.1} make it easy to demonstrate the
fact that fields of dimension $\le 1$ and without the property $C
_{1}$ are not uncommon for the class of algebraic extensions of any
local or global field.
\par
\medskip
\begin{rema}
\label{rema4.4} It is known that PAC fields have dimension $\le 1$
(see \cite{FJ}, Theorem~11.6.4, Corollary~11.2.5). Note also that
cd$(\mathcal{G}_{F}) \le 1$, for any field $F$ with dim$(F) \le
1$. Conversely, if $G$ is a profinite group with cd$(G) \le 1$,
then $G \cong \mathcal{G}_{F'}$, for some perfect PAC field $F
^{\prime }$ (see \cite{LvdD}, page~44; \cite{FJ},
Corollary~23.1.2). The question of whether almost perfect PAC
fields are of type $C _{1}$ seems to be open. As global fields are
almost perfect, it relates Problem \ref{prob} to the main topic of
this paper. A conjecture of Ax predicts, for perfect PAC fields,
that the answer to the question is affirmative. It has been proved
in characteristic zero (in \cite{Koll}), and for a perfect PAC
field $\Phi _{\ell }$ with $\mathcal{G}_{\Phi _{\ell }}$ a
pro-$\ell $-group, where $\ell \in \mathbb{P}$ (see \cite{W},
Sect. 3). These results imply one cannot prove that a field $F$
with dim$(F) \le 1$ and char$(F) = 0$ is not of type $C _{1}$,
using only topological invariants of $\mathcal{G}_{F}$.
\end{rema}
\par
\medskip
\begin{coro}
\label{coro4.5} With assumptions being as in the proof of Theorem
\ref{theo2.2}, for each finite extension $E ^{\prime }/E$ of odd
degree, there is a finite subset $m(E ^{\prime }) \subset
\mathbb{N}$, such that the forms $f _{n}$, $n \notin m(E ^{\prime
})$, are without nontrivial zeroes over $E ^{\prime }$.
\end{coro}
\par
\medskip
\begin{proof}
We first show that the norm group $N(\Lambda /E)$ equals $N(\Lambda E
^{\prime }/E ^{\prime }) \cap E ^{\ast }$, for any $\Lambda \in {\rm
Fe}(E)$ and each finite extension $E ^{\prime }$ of $E$ in an
algebraic closure of $E _{\rm sep}$, such that $\gcd (t, t') = 1$,
where $t = [\Lambda \colon E]$ and $t' = [E ^{\prime }\colon E]$. It
is easily verified that $\Lambda \cap E ^{\prime } = E$, $[\Lambda E
^{\prime }\colon E ^{\prime }] = [\Lambda \colon E]$ and the
restriction of the mapping $N _{\Lambda E'/E'}$ on $\Lambda $
coincides with $N _{\Lambda /E}$; therefore, $N(\Lambda /E) \subseteq
N(\Lambda E ^{\prime }/E ^{\prime }) \cap E ^{\ast }$. Note also that
$N _{\Lambda E'/E}$ equals the compositions $N _{\Lambda /E} \circ N
_{\Lambda E'/\Lambda }$ and $N _{E'/E} \circ N _{\Lambda E'/E'}$,
which yields $N(\Lambda E ^{\prime }/E) \subseteq N(\Lambda /E)$ and
$\alpha ^{t'} \in N(\Lambda E ^{\prime }/E)$, for all $\alpha \in
N(\Lambda E ^{\prime }/E ^{\prime }) \cap E ^{\ast }$. As $E ^{\ast
t} \subseteq N(\Lambda /E)$, these facts show that $N(\Lambda E
^{\prime }/E ^{\prime }) \cap E ^{\ast } \subseteq N(\Lambda /E)$, so
it turns out that $N(\Lambda /E) = N(\Lambda E ^{\prime }/E ^{\prime
}) \cap E ^{\ast }$, as claimed. It follows from this result that,
for any $n \in \mathbb{N}$, the form $f _{n}$ defined over $K$ as in
the proof of Theorem \ref{theo2.2}, can be defined essentially in the
same way over any finite extension $K _{n}$ of $K$ of degree
relatively prime to $p _{n}k _{n}(p _{n} - k _{n})$; hence, $f _{n}$
does not possess a nontrivial zero over $K _{n}$. Since, by (4.3),
$\gcd (p _{m}k _{m}(p _{m} - k _{m}), p _{n}k _{n}(p _{n} - k _{n}))$
is equal to $2 ^{m}$ if $m, n \in \mathbb{N}$ and $m < n$, this
proves Corollary \ref{coro4.5}.
\end{proof}
\par
\medskip
Our next result allows us to lift the restriction on $[E ^{\prime }\colon
E]$ in Corollary \ref{coro4.5}:
\par
\medskip
\begin{prop}
\label{prop4.6}
In the setting of Theorem \ref{theo2.2}, the algebraic extension $E$
of $K$ can be chosen so that {\rm dim}$(E) \le 1$ and there exist
$E$-forms $f _{n}$, $n \in \mathbb{N}$, without nontrivial zeroes over
$E$, which satisfy the following conditions:
\par
{\rm (a)} {\rm deg}$(f _{n}) = p _{n}$ and $f _{n}$ depends on $p
_{n}t _{n}$ variables, for each $n \in \mathbb{N}$ and some $p _{n},
t _{n} \in \mathbb{P}$, such that $t _{n} < p _{n}/3$;
\par
{\rm (b)} The sequence $p _{n}t _{n}$, $n \in \mathbb{N}$, increases
and consists of odd pairwise relatively prime numbers;
\par
{\rm (c)} For any finite extension $E ^{\prime }/E$, there is a
finite subset $m(E ^{\prime })$ of $\mathbb{N}$, such that none of the
forms $f _{n}$, $n \notin m(E ^{\prime })$, has a nontrivial zero
over $E ^{\prime }$.
\end{prop}
\par
\medskip
To prove Proposition \ref{prop4.6} we need the following lemma:
\par
\medskip
\begin{lemm}
\label{lemm4.7}
For each finite subset $P$ of $\mathbb{P}$, there exists $n(P) \in
\mathbb{N}$, such that every odd integer $N > n(P)$ is presentable as
a sum of three distinct prime numbers greater than any element of $P$.
\end{lemm}
\par
\medskip
\begin{proof}
This follows from the fact that for any $\alpha \in \mathbb{R}$ with
$0 < \alpha < 1$, there exists $M(\alpha ) \in \mathbb{R}$,
such that each integer $N > M(\alpha )$ equals the sum $N
= p _{1} + p _{2} + p _{3}$, for some $p _{i} \in \mathbb{P}$, $i = 1,
2, 3$, depending on $N$ so that $N ^{\alpha } < p _{1} < p _{2} < p
_{3}$. The fact itself has been established by Ax and deduced from
Vinogradov's theorem on the Ternary Goldbach Problem (see [4,
Lemma~2]).
\end{proof}
\par
\medskip\noindent
{\it Proof of Proposition \ref{prop4.6}}. Proceeding by induction on
$n$, and using Lemma \ref{lemm4.7}, one proves the existence of
$4$-tuples $(t _{n}, \theta _{n}, y _{n}, p _{n})$, $n \in
\mathbb{N}$, satisfying the following conditions, for each index
$n$:
\par
\medskip\noindent
(4.6) (a) $t _{n}, \theta _{n}, y _{n}$ and $p _{n}$ lie in $\mathbb
P$, and $t _{n} + \theta _{n} + y _{n} = p _{n}$;
\par
(b) $2 < t _{n} < \theta _{n} < y _{n}$ and $p _{n} < t _{n+1}$.
\par
\medskip
The rest of our proof goes along the lines drawn in the concluding
part of the proof of Theorem \ref{theo2.2} (after (4.3)), so we
present only its main steps and omit details. Put $S = \{t _{n},
\theta _{n}, y _{n}\colon n \in \mathbb{N}\}$ and $\Sigma = \{p
_{n}\colon n \in \mathbb{N}\}$. It follows from (4.6) that $S \cap
\Sigma = \emptyset $; also, by Lemma \ref{lemm4.1}, there is an
algebraic extension $E$ of $K$ with $S(E) = S$ and $\Sigma (E)
= \Sigma $. Therefore, dim$(E) \le 1$ and conditions (4.1) (b) hold,
which ensures the existence of cyclic extensions $T _{n}$, $\Theta
_{n}$ and $Y _{n}$, $n \in \mathbb{N}$, of $E$ in $E _{\rm ur}$ of
degrees $[T _{n}\colon E] = t _{n}$, $[\Theta _{n}\colon E] = \theta
_{n}$ and $[Y _{n}\colon E] = y _{n}$, for each $n$. Fix primitive
elements $\xi _{n} \in O _{v}(T _{n}) ^{\ast }$, $\eta _{n} \in O
_{v}(\Theta _{n}) ^{\ast }$, and $\delta _{n} \in O _{v}(Y _{n})
^{\ast }$ of $T _{n}/E$, $\Theta _{n}/E$ and $Y _{n}/E$,
respectively, so that the residue classes $\hat \xi _{n} \in \widehat
T _{n}$, $\hat \eta _{n} \in \widehat \Theta _{n}$ and $\hat \delta
_{n} \in \widehat Y _{n}$ be primitive elements of $\widehat T
_{n}/\widehat E$, $\widehat \Theta _{n}/\widehat E$ and $\widehat Y
_{n}/\widehat E$, respectively. Take a system $\overline X _{n} = X
_{n,i}$, $i = 1, \dots , t _{n}$, of algebraically independent
variables over $E$, and let $g _{n}(\overline X _{n})$ be the product
of the norms
\par\noindent
$$N _{T _{n}(\overline X _{n})/E(\overline X _{n})}(\sum _{i=1}
^{t_{n}} \xi _{n} ^{i-1}X _{n,i}), N _{\Theta _{n}(\overline X
_{n})/E(\overline X _{n})}(\sum _{i=1} ^{t _{n}} \eta _{n} ^{i-1}X
_{n,i}) $$ $${\rm and} \ N _{Y _{n}(\overline X _{n})/E(\overline X
_{n})}(\sum _{i=1} ^{t _{n}} \delta _{n} ^{i-1}X _{n,i}).$$
\par\noindent
Clearly, $g _{n}(\overline X _{n}) \in O _{v}(E)[\overline X _{n}]$
and the specializations $g _{n}(\overline \alpha _{n})$ of $g
_{n}(\overline X _{n})$, where $\overline \alpha _{n} = (\alpha
_{n,1}, \dots , \alpha _{n,t_{n}})$ and all $\alpha _{n,i} \in O
_{v}(E)$, are subject to the restrictions of (4.4) (here $k _{n}$ is
replaced by $t _{n}$). Also, $g _{n}(\overline X _{n})$
is an $E$-form of degree $p _{n}$, such that $v(g _{n}(\bar \beta
_{n})) \in p _{n}v(E)$ if $\bar \beta _{n} = (\beta _{n,1}, \dots ,
\beta _{n,t_{n}})$, $\beta _{n,i} \in E$, $i = 1, \dots , t _{n}$,
and $\beta _{n,i'} \neq 0$, for some $i'$.
Now fix $n$, observe that $v(E) \neq p _{n}v(E)$, take some $\pi _{n}
\in E ^{\ast }$ of value $v(\pi _{n}) > 0$ and $v(\pi _{n}) \notin p
_{n}v(E)$, and consider a system $\widetilde X = X _{n,i,j}\colon $
$i = 1, \dots , t _{n}$; $j = 1, \dots , p _{n}$, of $p _{n}t _{n}$
algebraically independent variables over $E$. It follows from the
preceding observations on $g _{n}(\overline X _{n})$ that the
polynomial $f _{n}(\widetilde X) = \sum _{j=1} ^{p _{n}} g
_{n}(\widetilde X _{n,j}).\pi _{n} ^{j}$ has the properties described
by (4.5) (with $t _{n}$ instead of $k _{n}$), where $\widetilde X
_{n,j} = X _{n,i,j}$, $i = 1, \dots , t _{n}$, for each fixed $j \in
\{1, \dots , p _{n}\}$. Moreover, arguing as in the proof of
Corollary \ref{coro4.5}, one obtains that $f _{n}$ is without
nontrivial zeroes over any extension of $E$ of degree relatively
prime to the product $p _{n}t _{n}\theta _{n}y _{n}$. Since (4.6)
implies $t _{n} < p _{n}/3$ and $\gcd (p _{m}t _{m}\theta _{m}y _{m},
p _{n}t _{n}\theta _{n}y _{n}) = 1$, for $m \in \mathbb{N}$, $m \neq
n$, this completes our proof.
\par
\medskip
\begin{rema}
\label{rema4.8} The problem of finding relations between Diophantine
properties of a field $E$ and the sequence cd$_{p}(\mathcal{G}_{E})$,
$p \in \mathbb{P}$, has attracted lasting interest in the study of
some modifications of the $C _{i}$-condition, introduced in
\cite{KK}. The research on this topic focuses on several conjectures
stated in \cite{KK}. One of them claims that a perfect field $E$ is
of type $C _{1}$, provided dim$(E) \le 1$ and $\mathcal{G}_{E}$ is a
pro-$p$-group, for some $p \in \mathbb{P}$ (and so suggests a
generalization of the Ax conjecture noted in Remark \ref{rema4.4}).
The stated generalization need not be true, see  \cite{CTM},
\cite{C-T}, but it seems to be unknown whether it holds in case $E$ is
an algebraic extension of $\mathbb{Q}$ or $\mathbb{Q} _{\ell }$, for
a fixed $\ell \in \mathbb{P}$. We refer the reader to \cite{W} and
\cite{I}, for more results on other modifications of the $C
_{i}$-condition, considered in \cite{KK}.
\end{rema}
\par
\medskip
Note finally that the present paper leaves open the question of
whether a Galois extension $E$ of a global or local field $K$ with
dim$(E) \le 1$ is a $C _{1}$-field. As shown in \cite{Ch2}, if $K$
is a local field, char$(\widehat K) = q$ and dim$(E) \le 1$, then
$S(E) \subseteq \{q\}$, i.e. $\widehat E(p) = \widehat E$, for all $p
\in \mathbb{P} \setminus \{q\}$. This, combined with Lemma
\ref{lemm3.5}, implies that if $v(E) \neq qv(E)$, then $I(E/K)$
contains a $K$-isomorphic copy of $K _{\rm ur}$. Since, by Lang's
theorem (referred to at the beginning of this Section), $K _{\rm ur}$
is of type $C _{1}$, it follows that so is $E$. On the other hand,
our research and the stated result of \cite{Ch2} indicate that the
method of proving Theorem \ref{theo2.2} does not lead to Galois
extensions $L$ of $K$ satisfying the conditions $v(L) = qv(L)$ and
dim$(L) \le 1$, which are not $C _{1}$-fields.

\par
\vskip0.4truecm\noindent \emph{Acknowledgements.} The author
wishes to thank the referee for a number of helpful suggestions,
used for improving the presentation of this research as a whole, and
particularly, of the proofs of Lemmas \ref{lemm3.4}, \ref{lemm3.5}
and Corollary \ref{coro4.5}. The research itself has partially been
supported by the Bulgarian National Science Fund under Grant KP-06 N
32/1 of 07.12.2019.

\end{document}